%% file: 596.tex
\documentstyle[psfig,diagram,12pt]{article}

\input myamssym.def
\input myamssym.tex

\input armacros

\catcode`\@=11

\def\@begintheorem#1#2{
\refstepcounter{subsection}
\rm \trivlist 
\item[\hskip \labelsep{\bf #1\ \thesubsection}]}

\def\@opargbegintheorem#1#2#3{
\refstepcounter{subsection}
\rm \trivlist
\item[\hskip \labelsep{\bf #1\ \thesubsection\  (#3):}]}

\catcode`\@=12

\title{Some independence results on reflection}

\author{James Cummings\\
        Hebrew University of Jerusalem\\
        {\tt cummings@math.huji.ac.il}
        \thanks{Supported by grant number 517/94 of the Israel Science Foundation.}
        \\
        Saharon Shelah\\
        Hebrew University of Jerusalem\\
        {\tt shelah@math.huji.ac.il}
        \thanks{
   Research supported by  the Israel Science Foundation,
   administered by The Israel Academy of Sciences and Humanities.
   Publication number 596.
                }
}

\begin{document}

\maketitle

\baselineskip=22pt

\begin{abstract}
  We prove that there is a certain degree of independence
 between stationary reflection phenomena at different cofinalities.
\end{abstract}

\section{Introduction}

 Recall that a stationary subset $S$ of a regular cardinal $\gk$
 is said to {\em reflect} at $\ga < \gk$ if $\cf(\ga) > \go$
 and $S \cap \ga$ is stationary in $\ga$. Stationary reflection
 phenomena have been extensively studied by set theorists,
 see for example \cite{Magidor}. 

\begin{definition} Let $\gk = \cf(\gk) < \gl = \cf(\gl)$.
$T^\gl_\gk =_{\rm def} \setof{\ga < \gl}{\cf(\ga)=\gk}$.
If $m < n < \go$ then $S^n_m =_{\rm def} \setof{\ga < \ha_n}{\cf(\ga) = \ha_m}$.
\end{definition}

 Baumgartner proved in \cite{JEB} that if $\gk$
 is weakly compact, GCH holds, and $\go < \gd = \cf(\gd) < \gk$ then
 forcing with the Levy collapse $Coll(\gd, <\gk)$ gives a model
 where for all $\gr < \gd$ and all stationary $T \subseteq T^\gk_\gr$
 the stationarity of $T$ reflects to some $\ga \in T^\gk_\gd$.
 In this last result all the cofinalities $\gr < \gd$ are on the same
 footing; we will build models where reflection holds for some
 cofinalities but fails badly for others.

 We introduce a more compact terminology for talking about reflection.

\begin{definition} Let $\gk = \cf(\gk) < \gl = \cf(\gl) < \gm = \cf(\gm)$.
\begin{enumerate}

\item  $Ref(\gm, \gl, \gk)$ holds iff for every stationary $S \subseteq T^\gm_\gk$
 there is a $\ga \in T^\gm_\gl$ with $S \cap \ga$ stationary in $\ga$.
 
\item  $Dnr(\gm, \gl, \gk)$ (Dense Non-Reflection) holds iff for every
 stationary  $S \subseteq T^\gm_\gk$ there is a stationary
$T \subseteq S$ such that for no $\ga \in T^\gm_\gl$ is
 $T \cap \ga$ stationary in $\ga$.

\end{enumerate}
\end{definition}

  We will use a variation on an idea from Dzamonja and Shelah's
  paper \cite{ShDz}.

\begin{definition}
 Let $\gk = \cf(\gk) < \gl = \cf(\gl) < \gm = \cf(\gm)$.
 $Snr(\gm, \gl, \gk)$ (Strong Non-Reflection) holds iff there is $F:T^\gm_\gk \lra \gl$
 such that for all $\ga \in T^\gm_\gl$ there is $C \subseteq \ga$
 closed and unbounded in $\ga$ with $F \restriction C \cap T^\gm_\gk$ strictly
 increasing.
\end{definition}
 
  As the name suggests, $Snr(\gm, \gl, \gk)$ is a strong
 failure of reflection. It is easy to see that if Jensen's Global
 $\square$ principle holds then $Snr(\gm, \gl, \gk)$ holds for all
 $\gk < \gl < \gm$; in some sense the strong non-reflection
 principle captures exactly that part of $\square$ which is useful for
 building non-reflecting stationary sets.

\begin{lemma} $Snr(\gm,\gl,\gk) \implies Dnr(\gm,\gl,\gk)$.
\end{lemma}
 
\begin{proof} Let $S \subseteq T^\gm_\gk$ be stationary,
 and let $F: T^\gm_\gk \lra \gl$ witness the strong non-reflection.
 Let $T \subseteq S$ be stationary such that $F \restriction T$
 is constant. Let $\ga \in T^\gm_\gl$ and let $C$ be a club in
 $\ga$ on which $F$ is strictly increasing, then $C$ meets $T$
 at most once and hence $T$ is non-stationary in $\ga$.
\end{proof}

 We will prove the following results in the course of this paper.

\smallskip

\noindent {\bf Theorem} \ref{thm1}: If the existence of a weakly compact cardinal is consistent, then
 $Ref(\ha_3, \ha_2, \ha_0) + Snr(\ha_3, \ha_2, \ha_1)$ is consistent.

\smallskip

\noindent {\bf Theorem} \ref{thm2}: If the existence of a measurable cardinal is consistent, then
  $Ref(\ha_3, \ha_2, \ha_0) + Snr(\ha_3, \ha_2, \ha_1) + Snr(\ha_3, \ha_1, \ha_0)$
is consistent.

\smallskip

\noindent {\bf Theorem} \ref{thm3}: If the existence of  a supercompact cardinal with a measurable
 above is consistent, then 
  $Ref(\ha_3, \ha_2, \ha_0) + Snr(\ha_3, \ha_2, \ha_1) + Ref(\ha_3, \ha_1, \ha_0)$
is consistent.

\smallskip

\noindent {\bf Theorem} \ref{zfc1}: If $Snr(\gm,\gl,\gk)$ and $\gk < \gk^* < \gl$ then
 $Snr(\gm, \gl, \gk^*)$.
 
\smallskip

\noindent {\bf Theorem} \ref{zfc2}: Let $\gk < \gl < \gm < \gn$. Then 
 $Ref(\gn,\gm,\gl) + Ref(\gn,\gl,\gk) \implies Ref(\gn,\gm,\gk)$
 and $Ref(\gn, \gm, \gk) + Ref(\gm, \gl, \gk) \implies Ref(\gn, \gl, \gk)$.

\smallskip

\noindent {\bf Theorem} \ref{thm4}:  If the existence of a weakly compact cardinal is consistent, then
 $Ref(\ha_3, \ha_2, \ha_1) + Dnr(\ha_3, \ha_2, \ha_0)$ is consistent.

\smallskip

  Theorems \ref{thm1}, \ref{thm2} and \ref{thm3} were proved by the first author. Theorems \ref{zfc1}
 and \ref{thm4} were proved
  by the second author, answering questions put to him by the first author.
  Theorem \ref{zfc2} was noticed by the first author (but has probably been observed many times).
  We would like to thank the anonymous referee for their very thorough reading of the first version
  of this paper.

\section{Preliminaries}

 We will use the idea of {\em strategic closure\/} of a partial
 ordering (introduced by Gray in \cite{Gray}).

\begin{definition} Let $\FP$ be a partial ordering, and let $\gee$
 be an ordinal.

\begin{enumerate}

\item
 The game $G(\FP, \gee)$ is played by two players
 I and II, who take turns to play elements $p_\ga$ of $\FP$
 for $ 0< \ga < \gee$, with player I playing at odd stages
 and player II at even stages (NB limit ordinals are even).

  The rules of the game are that the sequence that is played
 must be decreasing (not necessarily strictly decreasing),
 the first player who cannot make a move loses, and player II
 wins if play proceeds for $\gee$ stages.

\item
 $\FP$ is {\em $\gee$-strategically closed\/} iff 
 player II has a winning strategy in $G(\FP, \gee)$.

\item
 $\FP$ is {\em $<\gee$-strategically closed\/} iff for all $\gz < \gee$
 $\FP$ is  $\gz$-strategically closed.
 
\end{enumerate}

\end{definition}

  Strategic closure has some of the nice features of the standard
 notion of closure. For example a $(\gd+1)$-strategically closed
 partial ordering will add no $\gd$-sequences, and the property of
 being $(\gd + 1)$-strategically closed is preserved by forcing
 with $\le \gd$-support (see \cite{CuDzSh} for more information
 on this subject). We will need to know that under some circumstances
 we can preserve a stationary set by forcing with a poset that has
 a sufficient degree of strategic closure.

 The following is well-known.

\begin{lemma} Let GCH hold, let $\gl =\cf(\gl)$ and $\gk = \gl^+$.
Let $\gd = \cf(\gd) \le \gl$, and suppose that $\FP$ is
 $(\gd+1)$-strategically closed and $S$ is a stationary subset
 of $T^\gk_\gd$. Then $S$ is still stationary in $V^\FP$.
\end{lemma}

\begin{proof} Let $p \forces \mbox{``$\dot C$ is club in $\gk$''}$.
 Build
 $\seq{X_\ga: \ga < \gk}$
 a continuous increasing chain of elementary substructures
 of some large $H_\gth$ such that everything relevant is in $X_0$,
 $\card{X_\ga} = \gl$, ${}^{<\gl} X_\ga \subseteq X_{\ga+1}$.
 We make the remark here that this would not be possible for
 $\gl$ singular, and indeed the theorem can fail in that case
 (see \cite{ShSuccSing} for details).

  Now find some limit $\gg$ such that $X_\gg \cap \gk \in S$, clearly
 $\cf(\gg) = \gd$ and so ${}^{<\gd}X_\gg \subseteq X_\gg$.
 Let $\gb =_{\rm def} X_\gg \cap \gk$ and 
  fix $\seq{\gb_i: i < \gd}$ cofinal in $\gb$.
 Since $\FP, \gd \in X_0 \subseteq X_\gg$ we can find in $X_\gg$ a winning
 strategy $\gs$ for the game $G(\FP, \gd+1)$.
  
 Now we build a sequence $\seq{p_\ga: \ga \le \gd}$ such that 
\begin{enumerate}
\item For each even $\gb$, $p_\gb = \gs(\vec p \restriction \gb)$.
\item For $\gb < \gd$, $p_\gb \in \FP \cap X_\gg$.
\item For each $i < \gd$, there is $\gee$ such that
 $\gb_i < \gee < \gb$ and   $p_{2i+1} \forces \hat \gee \in \dot C$.
\end{enumerate}

  We can keep going because $X_\gg \prec H_\gth$,
 ${}^{<\gd} X_\gg \subseteq X_\gg$
 and $\gs$ is a winning strategy. At the end of the construction $p_\gd$
 is a refinement of $p_0$ which forces that $\gb$ is a limit point
 of $\dot C$, and we are done.
\end{proof} 

\section{Some consistency results}

  In this section we prove (starting from a weakly compact cardinal)
 the consistency of
 $ZFC + Ref(\ha_3, \ha_2, \ha_0) + Snr(\ha_3, \ha_2, \ha_1)$.
 We also show that together with this we can have either 
 $Ref(\ha_3, \ha_1, \ha_0)$ or $Snr(\ha_3, \ha_1, \ha_0)$. 

 We begin by defining a forcing ${\FP_{\rm Snr}}$ to enforce $Snr(\ha_3, \ha_2, \ha_1)$.

\begin{definition} $p$ is a condition in ${\FP_{\rm Snr}}$ iff $p$ is a function from
  a bounded subset of $S^3_1$ to $\ha_2$, and for every $\gg \in S^3_2$
 with $\gg \le \sup(\dom(p))$  there is $C$ club in $\gg$ such that
  $p \restriction C \cap S^3_1$ is strictly increasing.
 ${\FP_{\rm Snr}}$ is ordered by extension.
\end{definition}

  It is easy to see that ${\FP_{\rm Snr}}$ is $\go_2$-closed, and in fact
 that it is $\go_2$-directed closed. 

\begin{lemma} ${\FP_{\rm Snr}}$ is $(\go_2+1)$-strategically closed.
\end{lemma}

\begin{proof} We describe a winning strategy for player II
 in $G({\FP_{\rm Snr}}, \go_2+1)$. Suppose that $p_\ga$ is the condition
 played at move $\ga$. Let $\gb$ be an even ordinal, then
 at stage $\gb$ II will play as follows.

\smallskip

 Define $q_\gb =_{\rm def} \bigcup_{\ga < \gb} p_\ga$,
 $\gr_\gb =_{\rm def} \dom(q_\gb)$,
  and then let II play as follows: $p_\gb = q_\gb$ unless $\gb$ is limit and $\cf(\gb) = \ha_1$,
 in which case $p_\gb = q_\gb \cup \{ (\gr_\gb, \gb ) \}$. 

\smallskip

 The strategy succeeds because when play reaches stage $\go_2$,
  $\setof{\gr_\gb}{\gb < \go_2}$ is a club witnessing that
  $p_{\go_2}$ is a condition.  
\end{proof}

  This shows that ${\FP_{\rm Snr}}$ preserves cardinals and cofinalities
 up to $\ha_3$, from which it follows that
 $V^{\FP_{\rm Snr}} \models Snr(\ha_3,\ha_2, \ha_1)$. If GCH holds then
 $\card{{\FP_{\rm Snr}}} = \ha_3$, so ${\FP_{\rm Snr}}$ has the
 $\ha_4$-c.c.~and all
 cardinals are preserved.

    Now we define in $V^{\FP_{\rm Snr}}$ a forcing $\FQ$.
 This will enable us
 to embed ${\FP_{\rm Snr}}$ into the Levy collapse $Coll(\go_2, \go_3)$ in
 a particularly nice way. 

\begin{definition} In $V^{\FP_{\rm Snr}}$ let $F:S^3_1 \lra \ha_2$ be the
 function added by ${\FP_{\rm Snr}}$.  $q \in \FQ$ iff $q$ is a closed
 bounded subset of $\ha_3$, the order type of $q$ is less than
 $\ha_2$, and $F \restriction \lim(q) \cap S^3_1$ is strictly
 increasing. 
\end{definition}

   The aim of $\FQ$ is to add a club of order type $\go_2$ on
 which $F$ is increasing. It is clear that $\FQ$ is countably
 closed and collapses $\go_3$.

\begin{lemma} If GCH holds,  ${\FP_{\rm Snr}} * \dot \FQ$ is
 equivalent to $Coll(\ha_2, \ha_3)$.
\end{lemma} 

\begin{proof} 
 Since ${\FP_{\rm Snr}} * \dot \FQ$ has cardinality $\ha_3$ and collapses 
 $\ha_3$, it will suffice to show that it has an $\ha_2$-closed
 dense subset.
  To see this look at those conditions $(p, c)$ where $c \in V$,
 and $\max(c) = \sup(\dom(p))$. It is easy to see that this
 set is dense and $\ha_2$-closed.
\end{proof}

\begin{theorem} \label{thm1} Let $\gk$ be weakly compact, let GCH hold. 
 Define a two-step iteration by $\FP_0 =_{\rm def} Coll(\go_2, < \gk)$ and
  $\FP_1 =_{\rm def} (\FP_{\rm Snr})_{V^{\FP_0}}$. Then
 $V^{\FP_0 * \FP_1} \models   Ref(\ha_3, \ha_2, \ha_0)
                            + Snr(\ha_3, \ha_2, \ha_1)$.
\end{theorem}

\begin{proof} We will first give the proof for the case when $\gk$
 is measurable and then show how to modify it for the case when
 $\gk$ is just weakly compact.
 Assuming that $\gk$ is measurable, let $j: V \lra M$ be an
 elementary embedding into a transitive inner model with
 critical point $\gk$, where ${}^\gk M \subseteq M$.
 Notice that by elementarity and the closure of $M$,
 $j(\FP_0) = Coll(\go_2, <j(\gk))_M = Coll(\go_2, <j(\gk))_V$.

 Let $G$ be $\FP_0$-generic over $V$
 and let $H$ be $\FP_1$-generic over $V[G]$.
 We already know that $V[G*H] \models Snr(\ha_3, \ha_2, \ha_1)$,
 so let us assume that in $V[G*H]$ we have $S$ a stationary subset
 of $S^3_0$. To prove that $S$ reflects we will build a generic
 embedding with domain $V[G*H]$ extending $j$. Notice that since
 $\card{\FP_0 * \FP_1} = \gk$
 we can prove (by looking at canonical names) 
 that $V[G*H] \models {}^\gk M[G*H] \subseteq M[G*H]$, so in particular
 we have $S \in M[G*H]$.

 We start by forcing with $\FQ$ over $V[G*H]$ to get a generic
 object $I$. $H*I$ is generic over $V[G]$ for
 $(\FP*\FQ)_{V[G]}$ which is equivalent to $Coll(\go_2, \gk)$,
 so we can regard $G*H*I$ as being generic for $Coll(\go_2, \le \gk)$.
 Now let $J$ be $Coll(\go_2, [\gk, j(\gk)))$-generic over
 $V[G*H*I]$, then $G*H*I*J$ is $j(\FP_0)$-generic over $V$ (so
 a fortiori over $M$) and $j``G \subseteq G*H*I*J$ so that we
 can lift to get $j:V[G] \lra M[G*H*I*J]$.

 It remains to lift $j$ onto $V[G*H]$, for which we need to force
 a generic $K$ for $j(\FP_1)$ with the property that $j``H \subseteq K$.
 We will get $K$ by constructing a {\em master condition\/} in $j(\FP_0)$
 (that is, a condition refining all the conditions in $j``H$)
 and forcing below that master condition. A natural candidate for a master condition
 is $F =_{\rm def} \bigcup j``H$, where it  is easily seen
 (since $\crit(j)=\gk$ and $j \restriction \FP_1 = id$) that
 $F$ is the generic function from $\gk$ to $\ha_2$ added by
 $H$. The models $V[G]$ and $M[G*H*I*J]$ agree in their computations
 of $T^\gk_{\ha_1}$ and $T^\gk_{\ha_2}$, so $F$ is increasing on a club
 at all the relevant points below $\gk$.
 Now $\gk = (\ha_3)_{V[G]}$ is an ordinal of cofinality $\ha_2$
 in $M[G*H*I*J]$, but there is no problem here because $I$ has
 introduced a club in $\gk$ on which $F$ is increasing. Hence
 $F$ is a condition in $j(\FP_1)$ and we can force to get 
 $K \supseteq j``H$ as desired.

 We claim that $S$ is still stationary in $M[G*H*I*J*K]$.
 $I$ is generic for countably closed forcing, $J$ is generic
 for $\ha_2$-closed forcing and $K$ is generic for
 $\ha_2$-closed forcing so that $S$ (being a set of cofinality
 $\go$ ordinals) remains stationary. Now we argue as usual
 that since $j(S) \cap \gk = S$ and $\cf(\gk) = \go_2$
 in $M[G*H*I*J*K]$, there must exist $\ga < \gk$ in
 $V[G*H]$ such that $\cf(\ga) = \go_2$ and $S \cap \ga$
 is stationary in $\ga$.

 We promised at the start of this proof that we would show how to weaken the
 assumptions on $\gk$ from measurability to weak compactness. We will actually
 sketch two arguments, based on  two well-known characterisations of
 weak compactness. See Hauser's paper \cite{Kai} for detailed accounts
 of some similar arguments.

\begin{enumerate}

\item $\gk$ is weakly compact iff for every $A \subseteq V_\gk$ and every
 $\Pi^1_1$ formula $\phi$,
  $V_\gk \models \phi(A) \implies \exists \ga \; V_\ga \models \phi(A \cap V_\ga)$.

\item $\gk$ is weakly compact iff $\gk$ is strongly inaccessible and for every
 transitive $M$ such that $\card{M} = \gk$, ${}^{<\gk} M \subseteq M$ and
 $M$ models enough set theory there is a transitive set $N$ and an elementary
 embedding $j: M \lra N$ with $\crit(j) = \gk$.

\end{enumerate}

\medskip 

\noindent Argument 1: $\FP_0 * \FP_1 \subseteq V_\gk$, and so if $\dot S$
 is a name for a stationary subset of $\gk$ we can represent it by
 $S^* = \setof{(p, \ga)}{p \forces \hat \ga \in \dot S} \subseteq V_\gk$.
 The fact that $S$ is forced to be stationary can be written as a $\Pi^1_1$
 sentence (the universal second-order quantification is over names for clubs).
 Using the first characterisation given above we can find inaccessible $\ga < \gk$ such that
 $S^* \cap V_\ga$ is a $Coll(\go_2, < \ga) * (\FP_{Snr})_{V^{Coll(\go_2, <\ga)}}$
 name for a stationary set.  

%
%

  Now the argument is just like the one from a measurable, only with $\ga$
 playing the role of $\gk$ and $\gk$ replacing $j(\gk)$. We see that $S$
 has an initial segment $S \cap \ga$ which is stationary in a certain 
 intermediate generic extension, and just need to check that $S \cap \ga$
 remains stationary and that $\ga$ becomes a point of cofinality $\ha_2$.
 This is routine.

\medskip

\noindent Argument 2: Given a name $S$ for a stationary subset of $S^3_0$,
 build $S$ and $\FP_0 * \FP_1$ into an appropriate model $M$ of size
 $\gk$. Get $j$ as in the second characterisation as above.
 Repeat (mutatis mutandis) the argument from a measurable.
\end{proof}

   Having proved Theorem \ref{thm1}, it is natural to ask whether there
  is any connection between reflection to points of cofinality $\go_2$
  and reflection to points of cofinality $\go_1$. The following results provide
  a partial (negative) answer.

\begin{theorem} \label{thm2}
 Con($Ref(\ha_3,\ha_2,\ha_0) +
  Snr(\ha_3,\ha_2,\ha_1) + Snr(\ha_3,\ha_1,\ha_0)$) follows from the consistency
 of a measurable cardinal.
\end{theorem}

\begin{proof} Let $\gk$ be measurable. Without loss of generality 
 GCH holds (as we can move to the inner model $L[\gm]$). We will
 sketch a proof that we can force to get $Snr(\gk, \ha_1, \ha_0)$ without
 destroying the measurability of $\gk$. A more detailed argument
 for a very similar result is given in \cite{CuDzSh} (alternatively
 one can argue that because $\square$ holds in $L[\gm]$, $Snr(\gk, \ha_1, \ha_0)$
 is already true in that model). Let $j: V \lra M$ be the ultrapower
 map associated with some normal measure $U$ on $M$.

 We will do a reverse Easton iteration of length $\gk +2$, forcing at every
 regular cardinal $\ga \le \gk^+$ with $\FQ_\ga$, where $\FQ_\ga$ is the
 natural forcing to add a witness to $Snr(\ga, \ha_1, \ha_0)$ by initial segments.
 An easy induction shows that $\FQ_\ga$ is $<\ga$-strategically closed, the
 point being that the witnesses added below $\ga$ can be used to produce
 strategies in the game played on $\FQ_\ga$. The argument is exactly parallel to
 the proof of Lemma 6 in \cite{CuDzSh}.
 
 Let us break up the generic as $G * g * h$, where $G$ is $\FP_\gk$-generic,
 $g$ is $\FQ_\gk$-generic and $h$ is $\FQ_{\gk^+}$-generic.
 The key point is that $j(\FP_\gk)/G * g *h$ is $\gk^+$-strategically closed
 in $V[G * g *h]$, so that by GCH we can build $H \in V[G * g *h]$ which is
 $j(\FP_\gk)/G * g *h$-generic over $M[G * g * h]$.
 Since $j``G \subseteq G * g *h * H$, we can lift $j$ to get
 $j: V[G] \lra M[G * g * h * H]$. It is easy to see that $\bigcup j``g (= \bigcup g)$
 will serve as a master condition, so using GCH again we may find $g^+ \in V[G * g * h]$
 such that $j``g \subseteq g^+$ and $g^+$ is $j(\FP_\gk)$-generic over $M[G * g * h * H]$.
 Finally we claim that $j``h$ generates a $j(\FP_{\gk^+})$-generic filter over
 $M[G * g * h * H * g^+]$, because $\FP_{\gk^+}$ is distributive enough and
 $M[G * g * h * H * g^+] = \setof{j(F)(\gk)}{\dom(F) = \gk, F \in V[G*g]}$.
 Hence in $V[G * g * h]$ we can lift $j$ onto $V[G * g * h]$, so the
 measurability of $\gk$ is preserved.

 Now we just repeat the construction (from a measurable) of Theorem
 \ref{thm1}, and claim that in the final model $Snr(\ha_3,\ha_1,\ha_0)$ holds.
 The point is that if $F: T^\gk_{\ha_0} \lra \go_1$ witnesses the
 truth of $Snr(\gk, \ha_1, \ha_0)$, then in the final model $F$ witnesses
 $Snr(\ha_3, \ha_1, \ha_0)$ because the forcing from Theorem 1 does not
 change $T^\gk_{\ha_0}$ or $T^\gk_{\ha_1}$.
\end{proof}

 We can also go to the opposite extreme.

\begin{theorem} \label{thm3}
  Let $\gk$ be $\gl$-supercompact, where $\gl > \gk$ and $\gl$ is
 measurable. Let GCH hold. Then $Ref(\ha_3,\ha_2,\ha_0) + Snr(\ha_3,\ha_2,\ha_1)
 + Ref(\ha_3,\ha_1,\ha_0)$ holds
 in some forcing extension.

\end{theorem}

\begin{proof} We will start by forcing with $\FP =_{\rm def}Coll(\go_1, <\gk)$,
 after which $\gk$ is $\go_2$ and $\gl$ is still measurable. 
 Then we will do the construction of Theorem \ref{thm1}, that is we
 force with $(Coll(\gk, < \gl) * \FP_{\rm Snr})_{V^\FP}$. 
 Let $\FP_0 = Coll(\gk, < \gl)_{V^\FP}$
  and $\FP_1 = (\FP_{Snr})_{V^{\FP * \FP_0}}$.

 We need to check that $Ref(\ha_3,\ha_1,\ha_0)$ holds in the final model. To see
 this fix $j: V \lra M$ such that $\crit(j) = \gk$, $j(\gk) > \gl$
 and ${}^\gl M \subseteq M$. Let $G$, $H$ and $I$ be the generics
 for $\FP$, $\FP_0$ and $\FP_1$ respectively. 

 Since $\FP_0 * \FP_1$ is countably closed and has size $\gl$, we may
 find an embedding $i: \FP * \FP_0 * \FP_1 \lra j(\FP)$
 such that $i \restriction \FP = id$ and $j(\FP)/i``(\FP * \FP_0 * \FP_1)$
 is countably closed. Let $J$ be $j(\FP)/i``(\FP * \FP_0 * \FP_1)$-generic,
 then we can lift $j$ in the usual way to get
 $j: V[G] \lra M[G * H * I * J]$. Since $\FP_0 * \FP_1$ is a
 $\gk$-directed-closed forcing notion of size $\gl$, we may find a lower
 bound for $j`` (H*I)$ in $j(\FP_0 * \FP_1)$ and use it as a 
 master condition, forcing $K$ such that $j`` (H*I) \subseteq K$
 and lifting $j$ to $j: V[G * H * I] \lra M[G * H * I * J *K]$.

 Now suppose that in $V[G * H * I]$ we have $T \subseteq S^3_0$
 a stationary set. If $\gm = \bigcup j``\gl$ then we may
 argue as usual that $T \in M[G*H*I]$ and that
 $M[G*H*I] \models \hbox{$j(T) \cap \gm$ is stationary in $\gm$}$.
 Since $J * K$ is generic for countably closed forcing 
 and $j(T) \cap \gm \subseteq T^\gm_{\ha_0}$ this will still be true
 in $M[G * H * I * J *K]$,  and since $\cf(\gm) = \ha_1$
 in this last model we have by elementarity that $T$ reflects
 to some point in $S^3_1$.
\end{proof}

\section{Some ZFC results}

 In the light of the results from the last section it is natural to ask
 about the consistency of $Ref(\ha_3,\ha_2,\ha_1) + Snr(\ha_3,\ha_2,\ha_0)$. The first author
 showed by a rather indirect proof that this is impossible, and the
 second author observed that there is a simple reason for this.

\begin{theorem} \label{zfc1}
 If $Snr(\gm, \gl, \gk)$ and $\gk < \gk^* = \cf(\gk^*) < \gl$
 then $Snr(\gm, \gl, \gk^*)$. 
\end{theorem}

\begin{proof} Let $F: T^\gm_\gk \lra \gl$ witness $Snr(\gm, \gl, \gk)$.
 Define $F^*: T^\gm_{\gk^*} \lra \gl$  by
\[
    F^*: \gs \longmapsto \min \setof { \bigcup_{\ga \in C \cap T^\gm_\gk} F(\ga)}
                                     {\hbox{$C$ club in $\gs$}}.
\]

 We claim that $F^*$ witnesses $Snr(\gm, \gl, \gk^*)$. To see this, let $\gd \in T^\gm_\gl$
 and fix $D$ club in $\gd$ such that $F \restriction D \cap T^\gm_\gk$ is strictly
 increasing.
 Let $\gs \in \lim(D) \cap T^\gm_{\gk^*}$, and observe that
   $F^*(\gs) = \bigcup_{\ga \in D \cap T^\gm_\gk} F(\ga)$,
  because for any club $C$ in $\gs$ we have
\[
  \bigcup_{\ga \in C \cap T^\gm_\gk} F(\ga)
  \ge
  \bigcup_{\ga \in C \cap D \cap T^\gm_\gk} F(\ga)  
  =
  \bigcup_{\ga \in  D \cap T^\gm_\gk} F(\ga).
\]
  It follows  that $F^* \restriction \lim(D) \cap T^\gm_{\gk^*}$ is strictly
 increasing, because if $\gs_0, \gs_1 \in \lim(D) \cap T^\gm_{\gk^*}$ with
 $\gs_0 < \gs_1$ and $\gb$ is any point in $\lim(D) \cap T^\gm_\gk \cap (\gs_0, \gs_1)$
 then 
 $F^*(\gs_0) \le F(\gb) < F^*(\gs_1)$.
\end{proof}

  We also take the opportunity to record some other easy remarks,
 which put limits on the extent of the independence between different forms of
 reflection.

\begin{theorem} \label{zfc2}
 Let $\gk < \gl < \gm < \gn$ be regular. Then
\begin{enumerate}
\item $Ref(\gn, \gm, \gl) + Ref(\gn, \gl, \gk) \implies Ref(\gn, \gm, \gk)$.
\item $Ref(\gn, \gm, \gk) + Ref(\gm, \gl, \gk) \implies Ref(\gn, \gl, \gk)$.
\end{enumerate}
\end{theorem}

\begin{proof}
 For the first claim: let $S \subseteq T^\gn_\gk$ be stationary, and
 let us define  $T = \setof{\ga \in T^\gn_\gl}{\hbox{$S \cap \ga$ is stationary}}$.
 We claim that $T$ is stationary in $\gn$. To see this suppose $C$ is club in $\gn$
 and  disjoint from $T$, and consider $C \cap S$; this set must reflect at some
 $\gg \in T^\gn_\gl$, but then on the one hand $S \cap \gg$ is stationary (so
 $\gg$ is in $T$) while on the other hand $C$ is unbounded in $\gg$ (so $\gg \in C$),
 contradicting the assumption that $C$ and $T$ are disjoint.

 Now let $\gd \in T^\gn_\gm$ be such that $T \cap \gd$ is stationary in $\gd$.
 We claim that $S$ reflects at $\gd$. For if $D$ is club in $\gd$ then
 there is $\gg \in T \cap \lim(D)$ by the stationarity of $T \cap \gd$,
 and now since $D \cap \gg$ is club in $\gg$ and $S \cap \gg$ is stationary in
 $\gg$ there is $\gb \in D \cap S$.
 This proves the first claim. 

\smallskip

 For the second claim: $S \subseteq T^\gn_\gk$ be stationary, and let
 $\gd \in T^\gn_\gm$ be such that $S \cap \gd$ is stationary in
 $\gd$. Let $f: \gm \lra \gd$ be continuous increasing and cofinal in
 $\gd$, and let $S^* = \setof{\ga < \gm}{f(\ga) \in S}$. Then
 $S^*$ is stationary in $\gm$, and we can find $\gb \in T^\gm_\gl$
 such that $S^* \cap \gb$ is stationary in $\gb$. Now $f(\gb) \in T^\gn_\gl$
 and $S \cap f(\gb)$ is stationary in $f(\gb)$. This proves the
 second claim.
\end{proof}

\section{More consistency results}

 From the results in the previous section we saw in particular
 that we cannot have $Ref(\ha_3, \ha_2, \ha_1) + Snr(\ha_3, \ha_2, \ha_0)$.
 In this section we will see that $Ref(\ha_3, \ha_2, \ha_1) + Dnr(\ha_3, \ha_2, \ha_0)$
 is consistent.

 We will need some technical definitions
 and facts before we can start the main
 proof.

\begin{definition} Let $S$ be a stationary subset of $\ha_3$. We define
 notions of forcing $\FP(S)$, $\FQ(S)$
 and $\FR(S)$.
\begin{enumerate} 
\item $c$ is a  condition in $\FP(S)$ iff $c$ is a closed bounded subset of $\ha_3$
 such that $c \cap S = \emptyset$.
\item $d$ is a condition in $\FQ(S)$ iff $d$ is a function with $\dom(d) < \ha_3$,
 $d: \dom(d) \lra 2$, $d(\gg) = 1 \implies \gg \in S$, and
  for all $\ga \le \dom(d)$ if $\cf(\ga) > \go$ then there is $C \subseteq \ga$
 closed unbounded in $\ga$ such that $\gg \in C \implies d(\gg) = 0$.
\item $e$ is a condition in $\FR(S)$ iff $e$ is a closed bounded subset
 of $\ha_3$ such that for every point $\ga \in \lim(e)$  with $\cf(\ga) > \go$
 the set $S \cap \ga$ is non-stationary in $\ga$. 
\end{enumerate}
 In each case the conditions are ordered by end-extension.
\end{definition}

  The aims of these various forcings are respectively to kill the stationarity
  of $S$ ($\FP(S)$), to add a non-reflecting stationary subset of $S$
  ($\FQ(S)$), and to make $S$ non-reflecting on a closed unbounded set of
  points ($\FR(S)$). Notice that for some choices of $S$ the definitions of 
  $\FP(S)$ and $\FR(S)$ may not behave very well, for example if $S = \go_3$
  then $\FP(S)$ is empty and $\FR(S)$ only contains conditions of countable
  order type. Notice also that $\FQ(S)$ and $\FR(S)$ are countably
 closed.

\begin{lemma} \label{fred}
Let GCH hold. Let $S \subseteq S^3_0$ and suppose that there is
 a club $C$ of $\go_3$ such that $S \cap \ga$ is non-stationary 
 for all $\ga \in C$ with $\cf(\ga) > \go$. Let $T$ be a stationary subset 
 of $S^3_1$. Then $\FP(S)$ adds no $\go_2$-sequences of ordinals, and
 also preserves the stationarity of $T$.
\end{lemma}

\begin{proof} First we prove that $\FP(S)$ is $\go_2$-distributive.
 Let $\seq{D_\ga: \ga < \go_2}$ be a sequence of dense sets in $\FP(S)$,
 and let $c \in \FP(S)$ be a condition. Fix some large regular cardinal
 $\gth$ and let $c, C, \vec D, S \in X \prec H_\gth$ where $\card{X} = \go_2$,
 ${}^{\go_1} X \subseteq X$. Let $\gg$ be the ordinal $X \cap \go_3$,
 then $\cf(\gg) = \go_2$ by the closure of $X$. By elementarity
 it follows that $C$ is unbounded in $\gg$, so that $\gg \in C$
 and $S \cap \gg$ is non-stationary.
 
 Fix $B \subseteq \gg$
 closed unbounded in $\gg$ such that $B \cap S = \emptyset$ and
 $B$ has order type $\go_2$.
 Now we build a chain of conditions $c_\ga \in \FP(S) \cap X$ for $\ga < \go_2$
 such that $c_0 \le c$,  $c_{2\gb+1} \in D_\gb$ and $\max(c_{2\gb}) \in B$;
 we can continue at each limit stage because $B$ is disjoint from
 $S$ and $X$ is sufficiently closed. Finally we let
 $d =_{def} \bigcup_{\ga < \go_2} c_\ga \cup \{\gg\}$,
 then $d \le c$ and $d \in D_\gb$ for all $\gb < \go_2$. This
 shows that $\FP(S)$ is $\go_2$-distributive.

  The argument for the preservation of stationarity is similar.
 Let $\dot E$ be a $\FP(S)$-name for a closed unbounded subset
 of $\go_3$ and let $c \in \FP(S)$ be a condition.
 This time build $X$ such that $c, C, \dot E, S \in X \prec H_\gth$
 where $\card{X} = \go_1$, ${}^\go X \subseteq X$ and
 $\gd =_{def} \sup(X \cap \go_3) \in T$.
 Again $\gd \in C$, so $S \cap \gd$ is nonstationary.
 Choose $B \subseteq \gd$
 closed unbounded of order type $\go_1$ with $B \cap S = \emptyset$
 and build a chain of conditions $c_\ga \in \FP(S) \cap X$ such that $c_0 \le c$,
 $\max(c_{2\gb}) \in B$ and $c_{2\gb+1}$  forces that
 $\dot E \cap (\max(c_{2\gb}), \max(c_{2\gb+1})) \neq \emptyset$.
 Finally if $d =_{def} \bigcup_{\ga<\go_1} c_\ga \cup \{\gd\}$
 then $d \le c$ and $d \forces \hat \gd \in \hat T \cap \dot E$.
\end{proof}

  Now we describe a certain kind of forcing iteration.
It will transpire that all iterations of this type are
 $\go_4$-c.c.~and $(\go_2+1)$-strategically closed, so that
 in particular all cardinalities and cofinalities are preserved.

\begin{definition}
  Fix $F: \go_4 \times \go_3 \lra \go_4$ such that
 for all $\gb < \go_4$ the map $i \longmapsto F(\gb, i)$
 is a surjection from $\go_3$ onto $\gb$.

 $\FP_\gb$ is a {\em nice iteration\/} iff
 
\begin{enumerate}
\item $\gb \le \go_4$.
\item $\FP_\gb$ is an iteration of length $\gb$ with $\le \go_2$-supports.   
\item $\FQ_0 = \{ 0 \}$.
\item $\FQ_{2\gg+1}$ is $\FQ(\dot S_\gg)_{V^{P_{2\gg+1}}}$ where
 $\dot S_\gg$ is some $\FP_{2\gg+1}$-name for a stationary subset
 of $S^3_0$. Let $S^*_\gg$ be the non-reflecting stationary
 subset of $S_\gg$ which is added by $\FQ_{2\gg+1}$.
\item $\FQ_{2\gg}$ is $\FR(\dot R_\gg)$ where $R_\gg$ is the diagonal
 union of $\seq{S^*_{F(\gg, i)}: i < \go_3}$. That is
 $R_\gg = \setof{\gd \in S^3_0}{\exists i < \gd \; \gd \in S^*_{F(\gg, i)}}$.

\end{enumerate}

\end{definition}

  It is clear that an initial segment of a nice iteration is nice.
 Also every final segment of a nice iteration is countably closed, so that all the
 sets $S^*_\gg$ remain stationary throughout the iteration.
 The following remark will be useful later.

\begin{lemma} \label{keypoint}
 Let $\gg < \gd < \go_4$. 
 Then $R_\gg - R_\gd$ is non-stationary.
\end{lemma}

\begin{proof} Let $C$ be the closed and unbounded set of $i < \go_3$ such that 
\[
    \setof{F(\gg, j)}{ j < i} = \setof {F(\gd, j)}{ j < i} \cap \gg.
\]
 Let $i \in C \cap R_\gg$. Then for some $j < i$ we have
 $i \in S^*_{F(\gg, j)}$, and by the definition of $C$
 there is $k < i$ such that $F(\gg, j) = F(\gd, k)$.
 So $i \in S^*_{F(\gd, k)}$, $i \in R_\gd$, and we have
 proved that $C \cap R_\gg \subseteq R_\gd$.
\end{proof}

   We define a certain subset of $\FP_\gb$, which we call
  $\FP^*_\gb$.

\begin{definition}
 If $\FP_\gb$ is a nice iteration then $\FP^*_\gb$ is the set of conditions
 $p \in \FP_\gb$ such that
\begin{enumerate} 
\item $p(\gg) \in  \hat V$
 (that is, $p(\gg)$ is a canonical name for an object in $V$)
 for all $\gg \in \dom(p)$.
\item There is an ordinal $\gr(p)$ such that
\begin{enumerate}
\item $2\gd \in \dom(p) \implies \max(p(2\gd)) = \gr(p)$.
\item $2\gd+1 \in \dom(p) \implies \dom(p(2\gd+1)) = \gr(p)+1$.
\item $2\gd+1 \in \dom(p) \implies p(2\gd+1)(\gr(p)) = 0$.
\end{enumerate}
\item If $2\gd \in \dom(p)$ then $\forall i < \gr(p) \; 2F(\gd,i) +1 \in \dom(p)$.
\item If $2\gd +1 \in \dom(p)$ then $p \restriction 2\gd+1$ decides
 $S_\gd \cap \gr(p)$.
\end{enumerate}
\end{definition}

\begin{lemma} If  $p \in \FP^*_\gb$ and
 $2\gd \in \dom(p)$ then $p \restriction 2 \gd $ forces that $\gr(p) \notin R_\gd$.
\end{lemma}

\begin{proof} Let $\gr = \gr(p)$.
 By the definition of $\FP^*_\gb$, we see that
 $2 F(\gd, i) + 1 \in \dom(p)$ and
 $p(2 F(\gd, i) + 1)(\gr)=0$ for all $i < \gr$. This means that
 $p \forces  \gr \notin \dot S^*_{F(\gd, i)}$ for all $i < \gr$,
 which is precisely to say $p \restriction 2\gd \forces \gr \notin \dot R_\gd$. 
\end{proof}

\begin{lemma}
\label{contclos}
 Let $\FP_\gb$ be a nice iteration.
 Let $\gd \le \go_2$ be a limit ordinal and let $\seq{p_\gg: \gg < \gd}$
 be a decreasing sequence of conditions from $\FP^*_\gb$ 
 such that $\seq{\gr(p_\gg): \gg < \gd}$ is continuous and
 increasing.
  Define $q$ by setting $\dom(q) = \bigcup_{\gg < \gd} \dom(p_\gg)$,
  $\gr = \bigcup_{\gg < \gd} \gr(q_\gg)$, 
  $q(2 \gep) = \bigcup \setof{p_{\gg}(2 \gep)}{2 \gep \in \dom(p_\gg)} \cup \{ \gr \}$,
  $q(2 \gep + 1) = \bigcup \setof{p_{\gg}(2 \gep+1)}{2 \gep + 1 \in \dom(p_\gg)} \cup \{(\gr, 0) \}$.

   Then $q \in \FP^*_\gb$ and $\gr(q) = \gr$.
\end{lemma}

\begin{proof} Clearly it is enough to show that $q \in \FP_\gb$. Most of this
 is routine; the key points are that $q \restriction 2 \gep$
 forces that $R_\gep \cap \gr$ is non-stationary, and that
 $q \restriction 2 \gep + 1$ forces that $S^*_\gep \cap \gr$ is
 non-stationary.

 For the first point, observe  that for all
 sufficiently large $\gg < \gd$ we have $2 \gep \in \dom(p_\gg)$, so that
 by the last lemma
 $p_\gg \restriction 2 \gep \forces \gr(p_\gg) \notin R_\gep$; since
 $\seq{\gr(p_\gg): \gg < \gd}$ is continuous and $q \restriction 2 \gep$
 refines $p_\gg \restriction 2 \gep$, this implies that $q \restriction
 2 \gep$ forces that $R_\gep$ is not stationary in $\gr$.

 Similarly, $2 \gep + 1 \in \dom(p_\gg)$ for all large $\gg$, so that
  $p_\gg(2 \gep + 1)(\gr(p_\gg)) = 0$ for all large $\gg$. It follows
 immediately that $q \restriction (2 \gep + 1)$ forces that $S^*_\gep \cap \gr$
 is non-stationary.
\end{proof}

\begin{lemma}

 $\FP^*_\gb$ is dense in $\FP_\gb$, and $\FP_\gb$ is $(\go_2 +1)$-strategically
 closed.

\end{lemma}

\begin{proof} The proof is by induction on $\gb$. We prove first that
 $\FP^*_\gb$ is dense. 

\medskip

\noindent $\gb =0$: there is nothing to do.

\medskip

\noindent $\gb = 2\ga + 2$:
   fix $p \in \FP_\gb$. Since $\FP_{2\ga+1}$ is strategically closed
 and $\FP^*_{2\ga+1}$ is dense we may find $q_1 \in \FP^*_{2\ga+1}$
 such that $q_1 \le p \restriction (2\ga+1)$,
 $q_1$ decides $p(2\ga+1)$, and $\gr(q_1) > \dom(p(2\ga+1))$.
 
 Now we build a decreasing $\go$-sequence $\vec q$ of elements of $\FP^*_{2\ga+1}$ such that
 $\gr(q_n)$ is increasing and $q_{n+1}$ decides $\dot S_\ga \cap \gr(q_n)$; at each stage we use
 the strategic closure of $\FP_{2\ga+1}$ and the fact that $\FP^*_{2\ga+1}$ is a dense subset.
 After $\go$ steps we define $q$ as follows; let $\gr = \bigcup \gr(q_n)$,
 $\dom(q) = \bigcup_n \dom(q_n) \cup \{2\ga+1\}$, and 
\begin{enumerate}
\item For $\gg < \ga$, $q(2\gg+1) = \bigcup_n q_n(2\gg+1) \cup \{(\gr, 0)\}$.
\item For $\gg \le \ga$, $q(2\gg) = \bigcup_n q_n(2\gg) \cup \{\gr\}$.
\item $q(2\ga+1)(i) = p(2\ga+1)(i)$ if $i \in \dom (p(2\ga+1))$, and
 $0$ otherwise.
\end{enumerate}

  By the last lemma, $q \restriction 2 \ga +1 \in \FP^*_{2 \ga +1}$.
 It is routine to check that $q \in \FP^*_{2 \ga + 2}$
 and $\gr(q) = \gr$.

\medskip

\noindent $\gb = 2\ga+1$: this is exactly like the last case,
 except that now we demand 
$\forall i < \gr(q_n) \; 2 F(\ga, i) + 1 \in \dom(q_{n+1})$.

\medskip

\noindent $\gb$ is limit, $\cf(\gb) = \go_1$: Fix $\seq{\gb_i: i < \go_1}$
 which is continuous increasing and cofinal in $\gb$.
 Let $p \in \FP_\gb$. Find $q_0 \in \FP^*_{\gb_0}$ such that
 $q_0 \le p \restriction \gb_0$ and set
 $p_0 = q_0 \frown p \restriction [\gb_0, \gb)$.

 Now we define $q_i$ and $p_i$ by induction for $i \le \go_1$.
\begin{enumerate}
\item Choose $q_{i+1} \le p_i \restriction \gb_{i+1}$
 with $q_{i+1} \in \FP^*_{\gb_{i+1}}$, and then define
 $p_{i+1} = q_{i+1} \frown p \restriction [\gb_{i+1}, \gb)$.
\item For $i$ limit let $\gr_i = \bigcup_{j<i} \gr(q_j)$,
 $q_i(2\gg+1) = \bigcup_{j<i} q_j(2\gg+1) \cup \{ (\gr_i, 0) \}$,
 $q_i(2\gg) =  \bigcup_{j<i} q_j(2\gg) \cup \{\gr_i\}$.
 Then let $p_i = q_i \frown p \restriction [\gb_i, \gb_{i+1})$.
\end{enumerate}

  For $i < \go_1$ it is easy to see that $p_i$, $q_i$ are conditions.
 We claim that $q = q_{\go_1} \in \FP^*_\gb$. The only subtle point is to see
 that $q \in \FP_\gb$. Let $2 \gd +1 \in \dom(q)$. Then for all
 large $i$ we know $2 \gd +1 < \gb_i$, $2 \gd +1 \in \dom(q_i)$,
 so that in particular $q_i(2 \gd+1)(\gr_i) = 0$ for all large $i$.
 This means that $q(2 \gd + 1)$ is the characteristic function of a set
 which does not reflect at $\gr$, so is a legitimate condition in
 $\FQ_{2 \gd +1}$.  Similarly if $2 \gd \in \dom(q)$ then
 for all large $i$ we see that $q_i \restriction 2 \gd \forces \gr_i \notin R_\gd$,
 so that $q \restriction 2 \gd$ forces that the stationarity of $R_\gd$ does
 not reflect at $\gr$.

\medskip

\noindent $\gb$ is a limit, $\cf(\gb) = \go$ or $\go_2$:
 similar to the cofinality $\go_1$ case.

\medskip

\noindent $\cf(\gb) = \go_3$: easy because $\FP_\gb$ is the direct
 limit of the sequence $\seq{\FP_\gg: \gg < \gb}$.

  This concludes the proof that $\FP^*_\gb$ is dense. It is now easy to
 see that $\FP_\gb$ is $(\go_2 + 1)$-strategically closed; the strategy
 for player II is simply to play into the dense set $\FP^*_\gb$ at
 every successor stage, and to play a lower bound constructed as in
 Lemma \ref{contclos} at each limit stage.
\end{proof}

\begin{lemma}
  Let $\FP_{2 \gg}$ be a nice iteration of length less
 than $\go_4$. Then $\FP_{2 \gg} * \FP(\dot R_\gg)$
 is $(\go_2 + 1)$-strategically closed.
\end{lemma}

\begin{proof} This is just like the last lemma. 
\end{proof}

   Notice that the effect of forcing with $\FP(R_\gg)$ is to destroy the
 stationarity of all the sets $S^*_\gd$ for $\gd < \gg$.
 We are now ready to prove the main result of this section.

\begin{theorem} \label{thm4} If the existence of a weakly compact cardinal is consistent,
 then $Ref(\ha_3, \ha_2, \ha_1) + Dnr(\ha_3, \ha_2, \ha_0)$ is consistent.
\end{theorem}

\begin{proof} As in the proof of Theorem \ref{thm1}, we will first give a
 proof assuming the consistency of a measurable cardinal and then show
 how to weaken the assumption to the consistency of a weakly compact
 cardinal. We will need a form of ``diamond'' principle.

\begin{lemma} \label{diamond1} If $\gk$ is measurable and GCH holds, then in some forcing
 extension
\begin{enumerate}
\item  $\gk$ is  measurable.
\item  There exists a sequence $\seq{S_\ga: \ga < \gk}$ such that $S_\ga \subseteq \ga$ for all $\ga$,
 and for all $S \subseteq \gk$ there is a normal measure $U$ on $\gk$ such that if
 $j_U: V \lra M_U \simeq Ult(V, U)$ is the associated elementary embedding then
 $j_U(\vec S)_\gk = S$ (or equivalently, $\setof{\ga}{S \cap \ga = S_\ga} \in U$.
\end{enumerate}
\end{lemma}

\begin{proof}[Lemma \ref{diamond1}]  
   The proof is quite standard. For a similar construction given in more detail
   see \cite{Mitchell}.

  Fix $j: V \lra M$ the ultrapower map
  associated with some normal measure on $\gk$.
  Let $\FQ_\ga$ be the forcing whose conditions are sequences $\seq{T_\gb: \gb < \gg}$ where
  $\gg < \ga$ and $T_\gb \subseteq \gb$ for all $\gb$, ordered by end-extension. (This is
  really the same as the Cohen forcing $Add(\ga, 1)$). Let $\FP_{\gk+1}$ be a Reverse Easton
  iteration  of length
  $\gk +1$, where we force with $(\FQ_\ga)_{V^{\FP_\ga}}$ at each inaccessible $\ga \le \gk$.

  Let $G_\gk$ be $\FP_\gk$-generic over $V$ and let $g$  be $\FQ_\gk$-generic over $V[G_\gk]$.
  We will prove that the sequence given
  by $S_\ga = g(\ga)$ will work, by producing an appropriate $U$
  for each  $S \in V[G_\gk][g]$ with $S \subseteq \gk$. Let us fix such an $S$.

  By GCH and the fact that $j(\FP_\gk)/G_\gk * g$
  is $\gk^+$-closed in $V[G_\gk][g]$,
  we may build $H \in V[G_\gk][g]$ which is $j(\FP_\gk)/G_\gk * g$-generic over $M[G_\gk][g]$.
  Now for the key point: we define
  a condition $q \in \FQ_{j(\gk)}$ by setting $q(\ga) = g(\ga)$ for $\ga < \gk$ and
  $q(\gk) = S$. Then we build $h \ni q$ which is $\FQ_{j(\gk)}$-generic over
  $M[G_\gk][g][H]$, using GCH and the $\gk^+$-closure of $\FQ_{j(\gk)}$ in $V[G_\gk][g]$.
  
  To finish we define $j: V[G_\gk][g] \lra M[G_\gk][g][H][h]$ by $j: \dot \gt^{G_\gk * g} 
  \longmapsto j(\dot \gt)^{G_\gk * g * H * h}$, where this map is well-defined and
  elementary because $j``(G_\gk * g) \subseteq G_\gk * g * H * h$. The extended $j$
  is still an ultrapower by a normal measure ($U$ say) because
  $M[G_\gk][g][H][h] = \setof{j(F)(\gk)}{F \in V[G_\gk][g]}$. It is clear from the
  definition of this map $j$ that $j(g)(\gk) = h(\gk) = q(\gk) = S$, so  
  the model $V[G_\gk][g]$ is as required.

  This concludes the proof of Lemma \ref{diamond1}.
\end{proof}

   Fixing some reasonable coding of members of $H_{\ga^+}$ by subsets of $\ga$,
   we may write the diamond property in the following equivalent form: there
   is a sequence $\seq{x_\ga : \ga < \gk}$ such that for every $x \in H_{\gk^+}$
   there exists $U$ such that $j_U(\vec x)_\gk = x$. Henceforth we will assume
   that we have fixed a sequence $\vec x$ with this property.

   Now we describe a certain Reverse Easton forcing iteration of length $\gk + 1$. 
   It will be clear after the iteration is defined that it is $\ha_2$-strategically closed,
   so that in particular $\ha_1$ and $\ha_2$ are preserved.
   At stage $\ga < \gk$ we will force with $\FQ_\ga$, where $\FQ_\ga$ is trivial
   forcing unless 
\begin{enumerate}
\item $\ga$ is inaccessible.
\item $V^{\FP_\ga} \models \ga = \ha_3, \ga^+_V = \ha_4$
\item $x_\ga$ is a $\FP_\ga$-name for a nice iteration $\FP^\ga_{2\gd+1}$
 of some length $2 \gd + 1 < \ga^+$.
\end{enumerate}
   In this last case $\FQ_\ga$ is defined to be 
 $\FP^\ga_{2 \gd + 1} * \FP(R^\ga_\gd) * Coll(\ha_2, \ga)$.

   At stage $\gk$ (which will be $(\ha_3)_{V^{\FP_{\gk}}}$), we will
 do a nice iteration $\FQ_\gk$ of length $\gk^+$, with some book-keeping designed
 to guarantee that for every stationary $S \subseteq S^3_0$ in the final
 model there exists $T \subseteq S$ a non-reflecting stationary subset.
 So by design $Dnr(\ha_3, \ha_2, \ha_0)$ holds in the final model
 (and in fact so does $Dnr(\ha_3, \ha_1, \ha_0)$). It remains to be
 seen that $Ref(\ha_3, \ha_2, \ha_1)$ is true.

 Let $T \subseteq S^3_1$ be a stationary subset of $S^3_1$ in the final
 model. Since $\FQ_\gk$ has the $\gk^+$-c.c.~we may assume that
 $T$ is the generic extension by $\FP_\gk * (\FQ_\gk \restriction (2\gd + 1))$
 for some $\gd < \gk^+$. Let ${\overline \FQ} = \FQ_\gk \restriction (2\gd+1)$.
 Using the diamond property of $\vec x$ and the
 definition of the forcing iteration  we may find $U$ such that
\[
     j_U(\FP_\gk) = \FP_\gk * {\overline \FQ} * \FP(R_\gd) *  \FR_{\gk+1, j_U(\gk)}
\]
where $\FR_{\gk+1, j_U(\gk)}$ is the iteration above $\gk$. To save on notation, denote
 $j_U$ by $j$.

  Now if $M \simeq Ult(V, U)$ is the target model of $j$, then
  we may assume by the usual arguments  that $T \in M^{\FP_\gk * {\overline \FQ} }$.
  Notice that the last step in the iteration $\overline\FQ$ was
  to force with $\FR(R_\gd)$, that is to add a club of points at which the
  stationarity of $R_\gd$ fails to reflect. Applying Lemma \ref{fred}
  we see that $T$ is still stationary in the extension by 
  $\FP_\gk * {\overline \FQ} * \FP(R_\gd)$. Since $\FR_{\gk+1, j(\gk)}$ is
  $\ha_2$-strategically closed, $T$ will remain stationary in the extension by
  $j(\FP_\gk)$ (although of course $\gk$ will collapse to become some ordinal of
  cofinality $\ha_2$).

  To finish the proof we will build a generic embedding from 
  $V^{\FP_\gk * {\overline \FQ}}$ to
  $M^{j(\FP_\gk * {\overline \FQ})}$. It is easy to get
  $j: V^{\FP_\gk} \lra M^{j(\FP_\gk)}$, what is needed is a master
  condition for $\overline\FQ$ and $j$. Since $\gk$ is $\ha_3$ in $V^{\FP_\gk}$
  and $\FQ_\gk$ is an iteration with $\le \ha_2$-supports, it is clear what the
  condition should be, we just need to check that it works.

\begin{definition} Define $q$ by setting $\dom(q) = j``(2 \gd + 1)$,
 $q( j(2 \gg + 1)) = f_\gg \cup \{ (\gk, 0) \}$, $q(2 \gg) = C_\gg \cup \{ \gk \}$,
 where $f_\gg: \gk \lra 2$ is the function added by $\overline\FQ$
 at stage $2\gg+1$ and $C_\gg$ is the club added at stage $2 \gg$.
\end{definition}

  We claim that $q$ is a condition in $j({\overline \FQ})$.
 To see this we should first check that $f_\gg$ is the characteristic function
 of a non-stationary subset of $\gk$; this holds because at stage $\gk$ in
 the forcing $j(\FP_\gk)$ we forced with $\FP(R_\gd)$ and made $S^*_\gg$
 non-stationary for all $\gg < \gd$. We should also check that $R_\gg$
 is non-stationary in $\gg$, and again this is easy by Lemma \ref{keypoint}
 and the fact that $R_\gd$ has been made non-stationary.

  Forcing with $j({\overline \FQ})$
 adds no bounded subsets of $j(\gk)$, so that clearly $T$ is still stationary
 and $\cf(\gk)$ is still $\ha_2$ in the model $M^{j(\FP_\gk * {\overline \FQ})}$.
 By the familiar reflection argument, 
 there exists $\ga \in S^3_2$ such that $T \cap \ga$ is stationary in the
 model $V^{\FP_\gk * {\overline \FQ}}$. $T \cap \ga$ will
 still be stationary in $V^{\FP_\gk * \FQ_\gk}$, because the rest of the
 iteration $\FQ_\gk$ does not add any bounded subsets of $\ha_3$.
 We have proved that $Ref(\ha_3, \ha_2, \ha_1)$ holds in
 $V^{\FP_\gk * \FQ_\gk}$, which finishes the proof of
 Theorem \ref{thm4} using a measurable cardinal.

 It remains to be seen that we can replace the measurable cardinal by  a weakly
 compact cardinal. To do this we will use the following result of Jensen.

\begin{fact} Let $V = L$ and let $\gk$ be weakly compact. Then
 there exists a sequence $\seq{S_\ga: \ga < \gk}$ with $S_\ga \subseteq \ga$ for all
 $\ga$, such that for all $S \subseteq \gk$ and all $\gP^1_1$ formulae $\phi(X)$
 with one free second-order variable
\[
   V_\gk \models \phi(S) \implies (\exists \ga < \gk \; S_\ga = S \cap V_\ga, \,\,
 V_\ga \models \phi(S_\ga)).
\]
\end{fact}

   Using this fact we can argue exactly as in Argument 1 at the end of the
 proof of Theorem \ref{thm1}.
 This concludes the proof of Theorem \ref{thm4}.
\end{proof}

\end{document}

%% file: myamssym.tex
\expandafter\ifx\csname pre amssym.tex at\endcsname\relax \else  \fi
\expandafter\chardef\csname pre amssym.tex at\endcsname=\the\catcode`\@
\catcode`\@=11
\newsymbol\boxdot 1200
\newsymbol\boxplus 1201
\newsymbol\boxtimes 1202
\newsymbol\square 1003
\newsymbol\blacksquare 1004
\newsymbol\centerdot 1205
\newsymbol\lozenge 1006
\newsymbol\blacklozenge 1007
\newsymbol\circlearrowright 1308
\newsymbol\circlearrowleft 1309
\undefine\rightleftharpoons
\newsymbol\rightleftharpoons 130A
\newsymbol\leftrightharpoons 130B
\newsymbol\boxminus 120C
\newsymbol\Vdash 130D
\newsymbol\Vvdash 130E
\newsymbol\vDash 130F
\newsymbol\twoheadrightarrow 1310
\newsymbol\twoheadleftarrow 1311
\newsymbol\leftleftarrows 1312
\newsymbol\rightrightarrows 1313
\newsymbol\upuparrows 1314
\newsymbol\downdownarrows 1315
\newsymbol\upharpoonright 1316
 \let\restriction\upharpoonright
\newsymbol\downharpoonright 1317
\newsymbol\upharpoonleft 1318
\newsymbol\downharpoonleft 1319
\newsymbol\rightarrowtail 131A
\newsymbol\leftarrowtail 131B
\newsymbol\leftrightarrows 131C
\newsymbol\rightleftarrows 131D
\newsymbol\Lsh 131E
\newsymbol\Rsh 131F
\newsymbol\rightsquigarrow 1320
\newsymbol\leftrightsquigarrow 1321
\newsymbol\looparrowleft 1322
\newsymbol\looparrowright 1323
\newsymbol\circeq 1324
\newsymbol\succsim 1325
\newsymbol\gtrsim 1326
\newsymbol\gtrapprox 1327
\newsymbol\multimap 1328
\newsymbol\therefore 1329
\newsymbol\because 132A
\newsymbol\doteqdot 132B
 
\newsymbol\triangleq 132C
\newsymbol\precsim 132D
\newsymbol\lesssim 132E
\newsymbol\lessapprox 132F
\newsymbol\eqslantless 1330
\newsymbol\eqslantgtr 1331
\newsymbol\curlyeqprec 1332
\newsymbol\curlyeqsucc 1333
\newsymbol\preccurlyeq 1334
\newsymbol\leqq 1335
\newsymbol\leqslant 1336
\newsymbol\lessgtr 1337
\newsymbol\backprime 1038
\newsymbol\risingdotseq 133A
\newsymbol\fallingdotseq 133B
\newsymbol\succcurlyeq 133C
\newsymbol\geqq 133D
\newsymbol\geqslant 133E
\newsymbol\gtrless 133F
\newsymbol\sqsubset 1340
\newsymbol\sqsupset 1341
\newsymbol\vartriangleright 1342
\newsymbol\vartriangleleft 1343
\newsymbol\trianglerighteq 1344
\newsymbol\trianglelefteq 1345
\newsymbol\bigstar 1046
\newsymbol\between 1347
\newsymbol\blacktriangledown 1048
\newsymbol\blacktriangleright 1349
\newsymbol\blacktriangleleft 134A
\newsymbol\vartriangle 134D
\newsymbol\blacktriangle 104E
\newsymbol\triangledown 104F
\newsymbol\eqcirc 1350
\newsymbol\lesseqgtr 1351
\newsymbol\gtreqless 1352
\newsymbol\lesseqqgtr 1353
\newsymbol\gtreqqless 1354
\newsymbol\Rrightarrow 1356
\newsymbol\Lleftarrow 1357
\newsymbol\veebar 1259
\newsymbol\barwedge 125A
\newsymbol\doublebarwedge 125B
\undefine\angle
\newsymbol\angle 105C
\newsymbol\measuredangle 105D
\newsymbol\sphericalangle 105E
\newsymbol\varpropto 135F
\newsymbol\smallsmile 1360
\newsymbol\smallfrown 1361
\newsymbol\Subset 1362
\newsymbol\Supset 1363
\newsymbol\Cup 1264
 
\newsymbol\Cap 1265
 
\newsymbol\curlywedge 1266
\newsymbol\curlyvee 1267
\newsymbol\leftthreetimes 1268
\newsymbol\rightthreetimes 1269
\newsymbol\subseteqq 136A
\newsymbol\supseteqq 136B
\newsymbol\bumpeq 136C
\newsymbol\Bumpeq 136D
\newsymbol\lll 136E
 
\newsymbol\ggg 136F
 
\newsymbol\circledS 1073
\newsymbol\pitchfork 1374
\newsymbol\dotplus 1275
\newsymbol\backsim 1376
\newsymbol\backsimeq 1377
\newsymbol\complement 107B
\newsymbol\intercal 127C
\newsymbol\circledcirc 127D
\newsymbol\circledast 127E
\newsymbol\circleddash 127F
\newsymbol\lvertneqq 2300
\newsymbol\gvertneqq 2301
\newsymbol\nleq 2302
\newsymbol\ngeq 2303
\newsymbol\nless 2304
\newsymbol\ngtr 2305
\newsymbol\nprec 2306
\newsymbol\nsucc 2307
\newsymbol\lneqq 2308
\newsymbol\gneqq 2309
\newsymbol\nleqslant 230A
\newsymbol\ngeqslant 230B
\newsymbol\lneq 230C
\newsymbol\gneq 230D
\newsymbol\npreceq 230E
\newsymbol\nsucceq 230F
\newsymbol\precnsim 2310
\newsymbol\succnsim 2311
\newsymbol\lnsim 2312
\newsymbol\gnsim 2313
\newsymbol\nleqq 2314
\newsymbol\ngeqq 2315
\newsymbol\precneqq 2316
\newsymbol\succneqq 2317
\newsymbol\precnapprox 2318
\newsymbol\succnapprox 2319
\newsymbol\lnapprox 231A
\newsymbol\gnapprox 231B
\newsymbol\nsim 231C
\newsymbol\ncong 231D
\newsymbol\diagup 231E
\newsymbol\diagdown 231F
\newsymbol\varsubsetneq 2320
\newsymbol\varsupsetneq 2321
\newsymbol\nsubseteqq 2322
\newsymbol\nsupseteqq 2323
\newsymbol\subsetneqq 2324
\newsymbol\supsetneqq 2325
\newsymbol\varsubsetneqq 2326
\newsymbol\varsupsetneqq 2327
\newsymbol\subsetneq 2328
\newsymbol\supsetneq 2329
\newsymbol\nsubseteq 232A
\newsymbol\nsupseteq 232B
\newsymbol\nparallel 232C
\newsymbol\nmid 232D
\newsymbol\nshortmid 232E
\newsymbol\nshortparallel 232F
\newsymbol\nvdash 2330
\newsymbol\nVdash 2331
\newsymbol\nvDash 2332
\newsymbol\nVDash 2333
\newsymbol\ntrianglerighteq 2334
\newsymbol\ntrianglelefteq 2335
\newsymbol\ntriangleleft 2336
\newsymbol\ntriangleright 2337
\newsymbol\nleftarrow 2338
\newsymbol\nrightarrow 2339
\newsymbol\nLeftarrow 233A
\newsymbol\nRightarrow 233B
\newsymbol\nLeftrightarrow 233C
\newsymbol\nleftrightarrow 233D
\newsymbol\divideontimes 223E
\newsymbol\varnothing 203F
\newsymbol\nexists 2040
\newsymbol\Finv 2060
\newsymbol\Game 2061
\newsymbol\mho 2066
\newsymbol\eth 2067
\newsymbol\eqsim 2368
\newsymbol\beth 2069
\newsymbol\gimel 206A
\newsymbol\daleth 206B
\newsymbol\lessdot 236C
\newsymbol\gtrdot 236D
\newsymbol\ltimes 226E
\newsymbol\rtimes 226F
\newsymbol\shortmid 2370
\newsymbol\shortparallel 2371
\newsymbol\smallsetminus 2272
\newsymbol\thicksim 2373
\newsymbol\thickapprox 2374
\newsymbol\approxeq 2375
\newsymbol\succapprox 2376
\newsymbol\precapprox 2377
\newsymbol\curvearrowleft 2378
\newsymbol\curvearrowright 2379
\newsymbol\digamma 207A
\newsymbol\varkappa 207B
\newsymbol\Bbbk 207C
\newsymbol\hslash 207D
\undefine\hbar
\newsymbol\hbar 207E
\newsymbol\backepsilon 237F
\catcode`\@=\csname pre amssym.tex at\endcsname

%% file: armacros.tex
\newcommand{\ga}{\alpha}     
\newcommand{\gb}{\beta}      
\renewcommand{\gg}{\gamma}   
\newcommand{\gd}{\delta}     
\newcommand{\gep}{\epsilon}  
\newcommand{\gz}{\zeta}      
\newcommand{\gee}{\eta}      
\newcommand{\gth}{\theta}    
      
\newcommand{\gk}{\kappa}  
\newcommand{\gl}{\lambda}    
\newcommand{\gm}{\mu}        
\newcommand{\gn}{\nu}

\newcommand{\gr}{\rho}       
\newcommand{\gs}{\sigma}     
\newcommand{\gt}{\tau}

\newcommand{\go}{\omega}

\newcommand{\gP}{\Pi}

\newcommand{\ha}{\aleph}

\newcommand{\setof}[2]{{\{\; #1 \; \vert \; #2 \; \} } }
\newcommand{\seq}[1]{{\langle #1 \rangle} }
\newcommand{\card}[1]{{\vert #1 \vert} }

\newcommand{\forces}{\Vdash}

\renewcommand{\models}{\vDash}

\newcommand{\dom}{{\rm dom}}

\newcommand{\crit}{{\rm crit}}

\newcommand{\cf}{{\rm cf}}

\newcommand{\lra}{\longrightarrow}

\newtheorem{definition}{Definition}
\newtheorem{theorem}{Theorem}
\newenvironment{proof}{\noindent{\bf Proof:}}{\nopagebreak\mbox{}\newline\makebox[\textwidth]{\hfill$\blacklozenge$}\par\bigskip}

\newcommand{\implies}{\Longrightarrow}

\catcode`\@=11
\def\@begintheorem#1#2{\rm \trivlist \item[\hskip \labelsep{\bf #1\ #2:}]}
\def\@opargbegintheorem#1#2#3{\rm \trivlist
      \item[\hskip \labelsep{\bf #1\ #2\ (#3):}]}
\catcode`\@=12

\newtheorem{lemma}{Lemma}

\newtheorem{fact}{Fact}

\newcommand{\FP}{{\Bbb P}}
\newcommand{\FQ}{{\Bbb Q}}
\newcommand{\FR}{{\Bbb R}}

